\documentclass[10pt]{amsart}

\usepackage{amssymb}
\usepackage{amscd}
\usepackage{upref}

\theoremstyle{plain}
\newtheorem{theorem}{Theorem}[section]
\newtheorem{corollary}[theorem]{Corollary}
\newtheorem{lemma}[theorem]{Lemma}
\newtheorem{proposition}[theorem]{Proposition}

\newtheorem*{Theorem A}{Theorem A}
\newtheorem*{Theorem B}{Theorem B}
\newtheorem*{Theorem C}{Theorem C}

\theoremstyle{definition}
\newtheorem{definition}[theorem]{Definition}

\theoremstyle{remark}

\newtheorem{remark}[theorem]{Remark}

\numberwithin{equation}{section}

\begin{document}

\title[Cap products in string topology]{Cap products in string topology}

\author{Hirotaka Tamanoi }

\address{ Department of Mathematics,
 University of California Santa Cruz  \newline
\indent Santa Cruz, CA 95064}
\email[]{tamanoi@math.ucsc.edu}
\date{}
\begin{abstract}
Chas and Sullivan showed that the homology of the free loop space
$LM$ of an oriented closed smooth manifold $M$
admits the structure of a Batalin-Vilkovisky (BV) algebra equipped
with an associative product (loop product) and a Lie bracket
(loop bracket). We show that the cap product is compatible
with the above two products in the loop homology. Namely, the cap product with cohomology classes coming from $M$ via the circle action acts as derivations on the loop product as well as on the loop bracket. We show that Poisson identities and Jacobi identities hold for the cap product action, turning $H^*(M)\oplus \mathbb H_*(LM)$ into a BV algebra. Finally, we describe cap products in terms of the BV algebra structure in the loop homology.
\end{abstract}

\maketitle

\tableofcontents

\section{Introduction}

Let $M$ be a closed oriented smooth $d$-manifold. Let $D: H_*(M)
\xrightarrow{\cong} H^{d-*}(M)$ be the Poincar\'e duality map. Following a practice in string topology, we shift the homology grading downward by $d$ and let $\mathbb H_{-*}(M)=H_{d-*}(M)$. The Poincar\'e duality now takes the form $D:\mathbb H_{-*}(M) \xrightarrow{\cong} H^{*}(M)$. For a homology element $a$, let $|a|$ denote its 
$\mathbb H_*$-grading of $a$. 

The intersection product $\cdot$ in homology is defined as the Poincar\'e dual of the cup product. Namely, for $a,b\in \mathbb H_*(M)$, $D(a\cdot b)=D(a)\cup D(b)$. If $\alpha\in H^*(M)$ is dual to $a$, then $\alpha\cap b=a\cdot b$, and its Poincar\'e dual is $\alpha\cup D(b)$. Thus, through Poincar\'e duality, the intersection product, the cap product, and the cup product are all the same. In particular, the cap product and the intersection product commute:
\begin{equation}
\alpha\cap (b\cdot c)=(\alpha\cap b)\cdot c
=(-1)^{|\alpha||b|}b\cdot(\alpha\cap c). 
\end{equation}
In fact, the direct sum $H^*(M)\oplus \mathbb H_*(M)$ can be made into a graded commutative associative algebra with unit, given by $1\in H^0(M)$, using the cap and the cup product.

For an infinite dimensional manifold $N$, there is no longer
Poincar\'e duality, and geometric intersections of finite dimensional
cycles are all trivial. However, cap products can still be nontrivial
and the homology $H_*(N)$ is a module over the cohomology ring
$H^*(N)$.

When the infinite dimensional manifold $N$ is a free loop space $LM$
of continuous maps from the circle $S^1=\mathbb R/\mathbb Z$ to $M$,
the homology $\mathbb H_*(LM)=H_{*+d}(LM)$ has a great deal more structure. As before, $|a|$ denotes the $\mathbb H_*$-grading of a homology element $a$ of $LM$. Chas and
Sullivan \cite{CS} showed that $\mathbb H_*(LM)$ has a degree preserving associative graded
commutative product $\cdot$ called the loop product, a
Lie bracket $\{\ ,\ \}$ of degree $1$ called the loop bracket
compatible with the loop product, and the BV operator $\Delta$ of
degree $1$ coming from the homology $S^1$ action. These structures turn
$\mathbb H_*(LM)$ into a Batalin-Vilkovisky (BV) algebra. The purpose of this
paper is to clarify the interplay between the cap product with cohomology elements and the BV structure in $\mathbb H_*(LM)$.

Let $p: LM \rightarrow M$ be the base point map
$p(\gamma)=\gamma(0)$ for $\gamma\in LM$. For a cohomology class
$\alpha\in H^*(M)$ in the base manifold, its pull-back $p^*(\alpha)\in
H^*(LM)$ is also denoted by $\alpha$. Let $\Delta: S^1\times LM
\longrightarrow LM$ be the $S^1$-action map. This map induces a degree $1$ map $\Delta$ in homology given by $\Delta a=\Delta_*([S^1]\times a)$
for $a\in \mathbb H_*(LM)$. For a cohomology class $\beta\in H^*(LM)$, the formula 
$\Delta^*(\beta)=1\times\beta+\{S^1\}\times\Delta\beta$ defines a degree $-1$ map $\Delta$ in cohomology, where $\{S^1\}$ is the fundamental cohomology class. Although we use the same notation $\Delta$ in three different but closely related situations, what is meant by $\Delta$ should be clear in the context.  

\begin{Theorem A} Let $b,c\in \mathbb H_*(LM)$. 
The cap product with $\alpha\in H^*(M)$ graded commutes with the loop product. Namely 
\begin{equation}
\alpha\cap(b\cdot c)=(\alpha\cap b)\cdot
c=(-1)^{|\alpha||b|}b\cdot(\alpha\cap c).
\end{equation}

For $\alpha\in H^*(M)$, the cap product with $\Delta\alpha\in H^*(LM)$
acts as a derivation on the loop product and the
loop bracket\textup{:}
\begin{align}
(\Delta\alpha)\cap(b\cdot c)=(\Delta\alpha\cap b)\cdot c+
(-1)^{(|\alpha|-1)|b|}b\cdot (\Delta\alpha\cap c), \\
(\Delta\alpha)\cap\{b,c\}=\{\Delta\alpha\cap b,c\}+
(-1)^{|\alpha|-1)(|b|+1)}\{b,\Delta\alpha\cap c\}.
\end{align}

The operator $\Delta$ acts as a derivation on the cap product. Namely, for $\alpha\in H^*(M)$ and $b\in \mathbb H_*(LM)$.
\begin{equation}
\Delta(\alpha\cap b)=\Delta\alpha\cap b + (-1)^{|\alpha|} \alpha\cap
\Delta b.
\end{equation}
\end{Theorem A}
 
We recall that in the BV algebra $\mathbb H_*(LM)$, the following identities
are valid for $a,b,c\in \mathbb H_*(LM)$ \cite{CS}:
\begin{gather}
\Delta(a\cdot b)=(\Delta a)\cdot b + (-1)^{|a|}a\cdot \Delta b +
(-1)^{|a|}\{a,b\} 
\tag{BV identity} \\ 
\{a, b\cdot c\} =\{a,b\}\cdot c + (-1)^{|b|(|a|+1)}b\cdot\{a,c\} 
\tag{Poisson identity} \\ 
a\cdot b=(-1)^{|a||b|}b\cdot a,\qquad
\{a,b\}=-(-1)^{(|a|+1)(|b|+1)}\{b,a\}
\tag{Commutativity} \\
\{a,\{b,c\}\}=\{\{a,b\},c\}+(-1)^{(|a|+1)(|b|+1)}\{b,\{a,c\}\}
\tag{Jacobi identity}
\end{gather} 
Here, $\deg a\cdot b=|a|+|b|, \deg \Delta a=|a|+1$, and $\deg\{a,b\}=|a|+|b|+1$. 

We can extend the loop product and the loop bracket in $\mathbb H_*(LM)$ to
include $H^*(M)$ in the following way. For $\alpha\in H^*(M)$ and
$b\in \mathbb H_*(LM)$, we define the loop product and the loop bracket of 
$\alpha$ and $b$ by 
\begin{equation}
\alpha\cdot b=\alpha\cap b, \qquad
\{\alpha,b\}=(-1)^{|\alpha|}(\Delta\alpha)\cap b.
\end{equation} 

Furthermore, the BV structure in $\mathbb H_*(LM)$ can be extended to the direct sum $A_*=H^*(M)\oplus \mathbb H_*(LM)$ by defining the BV operator $\boldsymbol\Delta$ on $A_*$ to be trivial on $H^*(M)$ and to be the usual homological $S^1$ action $\Delta$ on $\mathbb H_*(LM)$. Here in $A_*$, elements in $H^k(M)$ are regarded as having homological degree $-k$. 

\begin{Theorem B} The direct sum $H^*(M)\oplus \mathbb H_*(LM)$ has a structure of a BV algebra. In particular, for $\alpha\in H^*(M)$ and $b,c\in \mathbb H_*(LM)$, the following form of
Poisson identity and the Jacobi identity hold\textup{:}
\begin{gather} 
\begin{split}
\{\alpha\cdot b,c\}&=\alpha\cdot\{b,c\}+
(-1)^{|b|(|c|+1)}\{\alpha,c\}\cdot b \\ &=\alpha\cdot\{b,c\} +
(-1)^{|\alpha||b|} b\cdot\{\alpha,c\},
\end{split} \\
\{\alpha,\{b,c\}\}=\{\{\alpha,b\},c\} +
(-1)^{(|\alpha|+1)(|b|+1)}\{b,\{\alpha,c\}\}.
\end{gather}
\end{Theorem B}

All the other possible forms of Poisson and Jacobi identities are
also valid, and the above two identities are the most nontrivial
ones. These identities are proved by using standard properties of the
cap product and the BV identity above in $\mathbb H_*(LM)$ relating the BV operator $\Delta$ and the loop bracket $\{\,,\,\}$, but without using Poisson identities nor Jacobi identities in the BV algebra $\mathbb H_*(LM)$.

The above identities may seem rather surprising, but they become
transparent once we prove the following result.

\begin{Theorem C} For $\alpha\in H^*(M)$, 
let $a=\alpha\cap [M]\in \mathbb H_*(M)$ be its Poincar\'e dual. Then for $b\in \mathbb H_*(LM)$,
\begin{equation}\label{cap}
\alpha\cap b=a\cdot b,\qquad (-1)^{|\alpha|}\Delta\alpha\cap
b=\{a,b\}.
\end{equation} 
More generally, for cohomology elements $\alpha_0,\alpha_1,\dotsc
\alpha_r\in H^*(M)$, let $a_0, a_1,\dotsc a_r\in \mathbb H_*(M)$ be their Poincar\'e duals. Then for $b\in \mathbb H_*(LM)$, we have
\begin{equation}\label{composition of derivations}
(\alpha_0\cup \Delta\alpha_1\cup\dotsm \cup \Delta\alpha_r)\cap b
=(-1)^{|a_1|+\dotsb+|a_r|}a_0\cdot\{a_1,\{a_2,\dotsc \{a_r,b\}\dotsb \}\}.
\end{equation}
\end{Theorem C}

Since the cohomology $H^*(M)$ and the homology $\mathbb H_*(M)$ are isomorphic through Poincar\'e duality and $\mathbb H_*(M)$ is a subring of $\mathbb H_*(LM)$, the first formula in \eqref{cap} is not surprising. However, the main difference between $H^*(M)$ and $\mathbb H_*(M)$ in our context is that the homology $S^1$ action $\Delta$ is trivial on $\mathbb H_*(M)\subset \mathbb H_*(LM)$, although cohomology $S^1$ action $\Delta$ is nontrivial on $H^*(M)$ and is related to loop bracket as in \eqref{cap}.  

Theorem A and Theorem C describes the cap product action of the cohomology $H^*(LM)$ on the BV algebra $\mathbb H_*(LM)$ for most elements in $H^*(LM)$. For example, for $\alpha\in H^*(M)$, the cap product with $\Delta\alpha$ is a derivation on the loop algebra $\mathbb H_*(LM)$ given by a loop bracket, and consequently the cap product with a cup product $\Delta\alpha_1\cup\cdots\cup\Delta\alpha_r$ acts on the loop algebra as a composition of derivations, which is equal to a composition of loop brackets, according to \eqref{composition of derivations}. If $H^*(LM)$ is generated by elements $\alpha$ and $\Delta\alpha$ for $\alpha\in H^*(M)$ (for example, this is the case when $H^*(M)$ is an exterior algebra, see Remark \ref{exterior algebra}), then Theorem C gives a complete description of the cap product with arbitrary elements in $H^*(LM)$ in terms of the BV algebra structure in $\mathbb H_*(LM)$. However, $H^*(LM)$ is general bigger than the subalgebra generated by $H^*(M)$ and $\Delta H^*(M)$. 

Since $\mathbb H_*(LM)$ is a BV algebra, in view of Theorem C, the validity of Theorem B may seem obvious. However, in the proof of Theorem B, we only used standard properties of the cap product and the BV identity. In fact, Theorem B gives an alternate elementary and purely homotopy theoretic proof of Poisson and Jacobi identities in $\mathbb H_*(LM)$, when at least one of the elements $a,b,c$ are in $\mathbb H_*(M)$. Similarly, Theorem C gives a purely homotopy theoretic interpretation of the loop product and the loop bracket if one of the elements are in $\mathbb H_*(M)$.

Our interest in cap products in string topology comes from an
intuitive geometric picture that cohomology classes in $LM$ are dual
to finite codimension submanifolds of $LM$ consisting of certain loop
configurations. We can consider configurations of loops
intersecting in particular ways (for example, two loops having their base points in common), or we can consider a family of loops
intersecting transversally with submanifolds of $M$ at
certain points of loops. In a given family of loops, taking the cap product with a cohomology class selects a subfamily of a certain loop configuration, which are ready for certain loop interactions. In this context, roughly speaking, composition of two interactions of loops correspond to the cup product of corresponding cohomology classes.

The organization of this paper is as follows. In section 2, we
describe a geometric problem of describing certain family of  intersection configuration of loops in terms of cap products. This gives a geometric motivation for the remainder of the paper. In section 3, we review the loop product in $\mathbb H_*(LM)$ in detail from the point of view of the intersection product in $\mathbb H_*(M)$. Here we pay careful attention to signs. In particular, we give a homotopy theoretic proof of graded commutativity in the BV algebra $\mathbb H_*(LM)$, which turned out to be not so trivial. In section 4, we prove compatibility
relations between the cap product and the BV algebra structure, and
prove Theorems A and B. In the last section, we prove Theorem C.

We thank the referee for numerous suggestions which lead to clarification and improvement of exposition.

\section{Cap products and intersections of loops}

Let $A_1,A_2,\dotsc A_r$ and $B_1,B_2,\dotsc B_s$ be oriented closed
submanifolds of $M^d$. Let $F\subset LM$ be a compact family of
loops. We consider the following question. 

\smallskip

\textbf{Question} : Fix $r$ points $0\le t_1^*,t_2^*,\dotsc,t_r^*\le
1$ in $S^1=\mathbb R/\mathbb Z$. Describe the homology class of the subset $I$ of the compact family $F$
consisting of loops $\gamma$ in $F$ such that $\gamma$ intersects
submanifolds $A_1,\dotsc A_r$ at time $t_1^*,\dotsc t_r^*$ and
intersects $B_1,\dotsc B_s$ at some unspecified time.

\smallskip

This subset $I\subset F$ can be described as follows. We consider the following diagram of an evaluation map and a projection map: 
\begin{equation} \label{eval}
\begin{CD}
\overset{s}{\overbrace{(S^1\times\dotsb\times S^1)}}\times LM @>{e}>>
  \overset{r}{\overbrace{M\times \dotsb\times M}}\times
\overset{s}{\overbrace{M\times\dotsb\times M}}  \\
@V{\pi_2}VV @.   \\
LM  @. 
\end{CD}
\end{equation}
given by $e\bigl((t_1,\dotsc t_s),\gamma\bigr)
=\bigl(\gamma(t_1^*),\dotsc \gamma(t_r^*), \gamma(t_1), \dotsc
\gamma(t_s)\bigr)$. Then the pull-back set
$e^{-1}(\prod_iA_i\times\prod_jB_j)$ is a closed 
subset of $S^1\times\dotsb\times S^1\times LM$. Let
\begin{equation*}
\tilde{I}=e^{-1}(\prod_iA_i\times \prod_jB_j)\cap 
(S^1\times\dotsb\times S^1\times F). 
\end{equation*} 
The set $I$ in
question is given by $I=\pi_2(\tilde{I})$. We want to understand this set $I$ homologically, including multiplicity. Although $e^{-1}(\prod_iA_i\times\prod_jB_j)$ is infinite dimensional, it has finite codimension in $(S^1)^r\times LM$. So we work cohomologically.

Let $\alpha_i, \beta_j\in H^*(M)$ be cohomology classes dual to $[A_i], [B_j]$ for $1\le i\le r$ and $1\le j\le s$. Then the subset $e^{-1}(\prod_iA_i\times \prod_jB_j)$ is dual to the cohomology class  $e^*(\prod_i\alpha_i\times\prod_j\beta_j)\in H^*((S^1)^s\times LM)$. Suppose the family $F$ is parametrized by a closed oriented manifold $K$ by an onto map $\lambda:K \longrightarrow F$ and let $b=\lambda_*([K])\in \mathbb H_*(LM)$ be the homology class of $F$ in $LM$. Then the homology class of $\tilde{I}$ in $(S^1)^s\times LM$ is given by 
\begin{equation}
[\tilde{I}]=e^*(\prod_i\alpha_i\times\prod_j\beta_j)\cap ([S^1\times\dotsb\times S^1]\times b). 
\end{equation} 
Note that the homology class $(\pi_2)_*([\tilde{I}])$ represents the homology class of $I$ with multiplicity. 

\begin{proposition}\label{loop intersection} With the above notation, $(\pi_2)_*([\tilde{I}])$ is given by the following formula in terms of the cap product or in terms of the BV structure\textup{:}
\begin{equation}
\begin{aligned}
(\pi_2)_*([\tilde{I}])&=(-1)^{\sum_jj|\beta_j|-s}
\bigl(\alpha_1\dotsm\alpha_r(\Delta\beta_1)\dotsm(\Delta\beta_s)\bigr)
\cap b \\
&=(-1)^{\sum_jj|\beta_j|-s}[A_1]\cdots[A_s]\cdot\{[B_1],\{\cdots\{[B_s],b\}\cdots\}
\in \mathbb H_*(LM). 
\end{aligned}
\end{equation}
\end{proposition} 
\begin{proof} The evaluation map $e$ in \eqref{eval} is 
given by the following composition.
\begin{multline*} 
\overset{s}{\overbrace{S^1\times\dotsb\times S^1}}\times LM 
\xrightarrow{1\times\phi} 
(S^1\times\dotsb\times S^1) \times
\overset{r+s}{\overbrace{LM\times\dotsb\times LM}} \\
\xrightarrow{T} 
\overset{r}{\overbrace{LM\times\dotsb\times LM}}\times
\overset{s}{\overbrace{(S^1\times LM)\times\dotsb\times(S^1\times LM)}}  \\
\xrightarrow{1^r\times\Delta^s}
(LM\times\dotsb\times LM) \times
(LM\times\dotsb\times LM)
\xrightarrow{p^r\times p^s} 
(M\times\dotsb\times M) \times
(M\times\dotsb\times M),
\end{multline*} 
where $\phi$ is a diagonal map, $T$ moves $S^1$ factors.  Since we
apply $(\pi_2)_*$ later, we only need terms in $e^*(\prod A_i\times
\prod B_j)$ containing the factor
$\{S^1\}\times\dotsb\times\{S^1\}$. Since
$\Delta^*p^*(\beta_j)=1\times p^*(\beta_j)
+\{S^1\}\times\Delta\beta_j$ for $1\le j\le s$, following the above
decomposition of $e$, we have
\begin{equation*}
e^*(\alpha_1\times\dotsb\times\alpha_r\times\beta_1\times\dotsb\times\beta_s)
=\varepsilon
\{S^1\}^s\times \bigl(\alpha_1\dotsm\alpha_r(\Delta\beta_1)\dotsm
(\Delta\beta_s)\bigr) + \text{ other terms},
\end{equation*}
where the sign $\varepsilon$ is given by $\varepsilon=
(-1)^{\sum_{\ell=1}^s(s-\ell)(|\beta_{\ell}|-1)
+s\sum_{\ell=1}^r|\alpha_{\ell}|}$.  Thus, taking the cap product with
$[S^1]^s\times b$ and applying $(\pi_2)_*$, we get
\begin{multline*}
{\pi_2}_*
\bigl(e^*(\alpha_1\times\dotsb\times\alpha_r\times
\beta_1\times\dotsb\times\beta_s)
\cap ([S^1]\times\dotsb\times[S^1]\times b)\bigr) \\
=(-1)^{\sum_{\ell=1}^{s}\ell|\beta_{\ell}|-s}
\alpha_1\dotsm\alpha_r(\Delta\beta_1)\dotsm(\Delta\beta_s)\cap b.
\end{multline*}
The second equality follows from the formula \eqref{composition of derivations}. 
\end{proof}

\begin{remark} In the diagram \eqref{eval}, 
in terms of cohomology transfer ${\pi_2}^!$ we have 
\begin{equation}
{\pi_2}^!
e^*(\alpha_1\times\dotsb\times\alpha_r\times\beta_1\times\dotsb\times\beta_s)
=\pm
\alpha_1\dotsm\alpha_r(\Delta\beta_1)\dotsm(\Delta\beta_s), 
\end{equation}
where ${\pi_2}^!(\alpha)\cap
b=(-1)^{s|\alpha|}{\pi_2}_*\bigl(\alpha\cap{\pi_2}_!(b)\bigr)$ for any
$\alpha\in H^*((S^1)^s\times LM)$ and $b\in \mathbb H_*(LM)$. Here
${\pi_2}_!(b)=[S^1]^s\times b$.
\end{remark}


\section{The intersection product and the loop product}

Let $M$ be a closed oriented smooth $d$-manifold.  The loop product in
$\mathbb H_*(LM)$ was discovered by Chas and Sullivan \cite{CS}, in terms of
transversal chains. Later, Cohen and Jones \cite{CJ} gave a homotopy
theoretic description of the loop product. The loop product is a
hybrid of the intersection product in $\mathbb H_*(M)$ and the Pontrjagin
product in the homology of the based loop spaces $H_*(\Omega M)$. In
this section, we review and prove some properties of the loop product
in preparation for the next section. Our treatment of the loop product
follows \cite{CJ}. However, we will be precise with signs and give a
homotopy theoretic proof of the graded commutativity of the loop
product, which \cite{CJ} did not include. For the Frobenius compatibility formula with careful discussion of signs, see \cite{T2}. For homotopy theoretic deduction of the BV identity, see \cite{T3}. 

For our purpose, the free loop space $LM$ is the space of {\it continuous} maps from $S^1=\mathbb R/\mathbb Z$ to $M$. Our discussion is homotopy theoretic and does not require smoothness of loops, although we do need smoothness of $M$ which is enough to allows us to have tubular neighborhoods for certain submanifolds in the space of continuous loops. Recall that the space $LM$ of continuous loops can be given a structure of a smooth manifold. See the discussion before Definition \ref{definition of loop product}. 

Let $p: LM \longrightarrow M$ be the base
point map given by $p(\gamma)=\gamma(0)$. Let $s: M \longrightarrow
LM$ be the constant loop map given by $s(x)=c_x$, where $c_x$ is the
constant loop at $x\in M$. Since $p_*\circ s_*=1$, $\mathbb H_*(M)$ is
contained in $\mathbb H_*(LM)$ through $s_*$ and we often regard $\mathbb H_*(M)$ as a subset of $\mathbb H_*(LM)$.

We start with a discussion on the intersection ring $\mathbb H_*(M)$ and later we compare it with the loop homology algebra $\mathbb H_*(LM)$. An exposition on intersection products in homology of manifolds can be found on Dold's book \cite{D}, Chapter VIII, \S13. Our sign convention (which follows Milnor \cite{M}) is slightly different from Dold's. We give a fairly detailed discussion of the intersection ring $\mathbb H_*(M)$ because the
discussion for the loop homology algebra goes almost in parallel, and
because our choice of the sign for the loop product comes from and is
compatible with the intersection product in $\mathbb H_*(M)$. Compare formulas \eqref{intersection product} and \eqref{loop product}.

Those who are familiar with intersection product and loop products can skip this section after checking Definition \ref{definition of loop product}. 

Let $D: \mathbb H_*(M) \xrightarrow{\cong} H^{d-*}(M)$ be the Poincar\'e duality map such that $D(a)\cap[M]=a$ for $a\in \mathbb H_*(M)$.  We discuss two ways to define
intersection product in $\mathbb H_*(M)$. The first method is the official one
and we simply define the intersection product as the Poincar\'e dual
of the cohomology cup product. Thus, $D(a\cdot b)=D(a)\cup D(b)$ for
$a,b\in \mathbb H_*(M)$. For example, we have $a\cdot
b=(-1)^{|a||b|}b\cdot a$.

The second method uses the transfer map induced from the diagonal map
$\phi: M \longrightarrow M\times M$. Let $\nu$ be the normal bundle to
$\phi(M)$ in $M\times M$, and we orient $\nu$ by $\nu\oplus
\phi_*(TM)\cong T(M\times M)|_{\phi(M)}$. Let $u'\in
H^d(\phi(M)^{\nu})$ be the Thom class of $\nu$. Let $N$ be a closed
tubular neighborhood of $\phi(M)$ in $M\times M$ so that $D(\nu)\cong
N$, where $D(\nu)$ is the associated closed disc bundle of $\nu$. Let
$\pi: N \longrightarrow M$ be the projection map. Then the above Thom
class can be thought of as $u'\in \tilde{H}^d(N/\partial N)$, and we
have the following commutative diagram, where $c: M\times M
\longrightarrow N/\partial N$ is the Thom collapse map, and $\iota_N$
and $j$ are obvious maps.
\begin{equation}
\begin{CD} 
H^d(N, N-\phi(M)) @>{\cong}>> H^d(N,\partial N)\ni u' \\
@A{\cong}A{\iota_N^*}A   @VV{c^*}V   \\
u''\in H^d\bigl(M\times M, M\times M-\phi(M)\bigr) @>{j^*}>>  H^d(M\times M)\ni u
\end{CD}
\end{equation} 
Let $u''\in H^d\bigl(M\times M, M\times M-\phi(M)\bigr)$ and $u\in
H^d(M\times M)$ be the classes corresponding to the Thom class. We
have $u=c^*(u')=j^*(u'')$. This class $u$ is characterized by the
property $u\cap[M\times M]=\phi_*([M])$, and $\phi^*(u)=e_M\in H^d(M)$
is the Euler class of $M$. See for example section 11 of \cite{M}. The
transfer map $\phi_!$ is defined as the following composition:
\begin{equation}
\phi_!: H_*(M\times M) \xrightarrow{c_*} \tilde{H}_*(N/\partial N)
\xrightarrow[\cong]{u'\cap(\ )} H_{*-d}(N) \xrightarrow[\cong]{\pi_*}
H_{*-d}(M).
\end{equation} 

For a homology element $a$, let $|a|'$ denote its regular homology degree of $a$, so that we have $a\in H_{|a|}(M)$ and $|a|'=|a|+d$. 

\begin{proposition}\label{properties of phi} Suppose $M$ is a connected 
oriented closed $d$-manifold with a base point $x_0$. The transfer map
$\phi_!: H_*(M\times M) \longrightarrow H_{*-d}(M)$ satisfies the
following properties. For $a,b\in H_*(M)$,
\begin{align}
\phi_*\phi_!(a\times b)&=u\cap(a\times b) \\ \phi_!\phi_*(a\times
b)&=\chi(M)[x_0]
\end{align}
For $\alpha\in H^*(M\times M)$ and $b,c\in H_*(M)$, we have 
\begin{equation}
\phi_!\bigl(\alpha\cap(b\times c)\bigr)
=(-1)^{d|\alpha|}\phi^*(\alpha)\cap\phi_!(b\times c).
\end{equation} 
The intersection product and the transfer map coincide up to a sign. 
\begin{equation}\label{intersection product}
a\cdot b=(-1)^{d(|a|'-d)}\phi_!(a\times b).
\end{equation}
\end{proposition} 
\begin{proof} For the first identity, we consider the following 
commutative diagram, where $M^2$ denotes $M\times M$. 
\begin{equation*}
\begin{CD}
H_*(M^2) @>{c_*}>> H_*(N,\partial N) @>{u'\cap (\ )}>{\cong}> H_{*-d}(N) 
@>{\pi_*}>{\cong}>  H_{*-d}(M)  \\
@|   @V{\cong}V{{\iota_N}_*}V  @VV{{\iota_N}_*}V   @VV{\phi_*}V   \\
H_*(M^2)  @>{j_*}>> H_*\bigl(M^2, M^2-\phi(M)\bigr) 
@>{u''\cap(\ )}>{\cong}>  H_{*-d}(M^2) @=  H_{*-d}(M^2)
\end{CD}\end{equation*}
The commutativity implies that for $a,b\in H_*(M)$, we have $\phi_*\phi_!
=u''\cap j_*(a\times b)=j^*(u'')\cap(a\times b)=u\cap(a\times b)$. 

To check the second formula, we first compute
$\phi_*\phi_!\phi_*([M])$. By the first formula,
$\phi_*\phi_!\phi_*([M])=u\cap\phi_*([M])=\phi_*(\phi^*(u)\cap[M])$. Since
$\phi^*(u)$ is the Euler class $e_M$, this is equal to
$\phi_*(e_M\cap[M])=\chi(M)[(x_0,x_0)]$. Since $M$ is assumed to be
connected, $\phi_*$ is an isomorphism in $H_0$. Hence
$\phi_!\phi_*([M])=\chi(M)[x_0]\in H_0(M)$.

For the next formula, we examine the following commutative diagram. 
\begin{equation*}
\begin{CD} 
H_*(M^2) @>{c_*}>> H_*(N,\partial N) @>{u'\cap(\ )}>{\cong}> H_{*-d}(N)
@<{\iota_*'}<{\cong}< H_{*-d}(M) \\ 
@V{\alpha\cap(\ )}VV
@V{\iota_N^*(\alpha)\cap(\ )}VV @V{\iota_N^*(\alpha)\cap(\ )}VV
@V{\iota^*(\alpha)}VV \\ 
H_{*-|\alpha|}(M^2) @>{c_*}>>
H_{*-|\alpha|}(N,\partial N) @>{u'\cap(\ )}>{\cong}> H_{*-d-|\alpha|}(N)
@<{\iota_*'}<{\cong}< H_{*-d-|\alpha|}(M)
\end{CD}
\end{equation*}
where $\iota': M\rightarrow N$ is an inclusion map and
$\iota_*'=(\pi_*)^{-1}$. The middle square commutes up to
$(-1)^{|\alpha|d}$. The commutativity of this diagram immediately
implies that $\iota^*(\alpha)\cap\phi_!(a\times b)
=(-1)^{|\alpha|d}\phi_!\bigl(\alpha\cap(a\times b)\bigr)$.

For the last identity, we apply $\phi_*$ on both sides and
compare. Since $a\cdot b=\phi^*(D(a)\times D(b))\cap[M]$, we have
\begin{align*} 
\phi_*(a\cdot b)&=\bigl(D(a)\times D(b)\bigr)\cap \phi_*([M]) \\
&=\bigl(D(a)\times D(b)\bigr)\cap \bigl(u\cap[M]\bigr) \\
&=(-1)^{d(|a|'-d)}u\cap(a\times b)=(-1)^{d(|a|'-d)}\phi_*\phi_!(a\times b).
\end{align*} 
Since $\phi_*$ is injective, we have $a\cdot
b=(-1)^{d(|a|'-d)}\phi_!(a\times b)$.
\end{proof} 

These two intersection products, one defined using the Poincar\'e
duality, and the other using Pontrjagin Thom construction, differ only
in signs. However, the formulas for graded commutativity take
different forms.
\begin{align}
a\cdot b&=(-1)^{(d-|a|')(d-|b|')}b\cdot a \\
\phi_!(a\times b)&=(-1)^{|a|'|b|'+d}\phi_!(b\times a)
\end{align}
The sign $(-1)^d$ in the second formula above comes from the fact that
the Thom class $u\in H^d(M\times M)$ satisfies $T^*(u)=(-1)^du$, where
$T$ is the switching map of factors.

Next we turn to the loop product in $H_*(LM)$. We consider the
following diagram.
\begin{equation}
\begin{CD}
LM\times LM @<{j}<< LM\times_M LM  @>{\iota}>> LM  \\
@V{p\times p}VV    @V{q}VV        @.  \\
M\times M @<{\phi}<<  M @.  
\end{CD}
\end{equation}
where $LM\times_M LM=(p\times p)^{-1}(\phi(M))$ consists of pairs of
loops $(\gamma,\eta)$ with the same base point, and
$\iota(\gamma,\eta)=\gamma\cdot\eta$ is the product of composable
loops. Let $\widetilde{N}=(p\times p)^{-1}(N)$ and let $\tilde{c}:
LM\times LM \longrightarrow \widetilde{N}/\partial\widetilde{N}$ be the Thom collapse map. Let $\tilde{\pi}: \widetilde{N} \longrightarrow
LM\times_MLM$ be a projection map defined as follows. For
$(\gamma,\eta)\in\widetilde{N}$, let their base points be $(x,y)\in
N$. Let $\pi(x,y)=(z,z)\in \phi(M)$. Since $N\cong D(\nu)$ 
has a bundle structure,
let $\ell(t)=(\ell_1(t),\ell_2(t))$ be the straight ray in the fiber
over $(z,z)$ from $(z,z)$ to $(x,y)$. Then let
$\tilde{\pi}\bigl((\gamma,\eta)\bigr)=(\ell_1\cdot\gamma\cdot\ell_1^{-1},
\ell_2\cdot\eta\cdot\ell_2^{-1})$. By considering $\ell_{[t,1]}$, we
see that $\tilde{\pi}$ is a deformation retraction. 

In fact, more is true. Stacey (\cite{St}, Proposition 5.3) showed that when $L_{\text{smooth}}M$ is the space of {\it smooth} loops, $\widetilde{N}$ has an actual structure of a tubular neighborhood of $LM\times_MLM$ inside of $LM\times LM$ equipped with a diffeomorphism  $p^*\bigl(D(\nu)\bigr)\cong \widetilde{N}$. His proof only uses the smoothness of $M$ and exactly the same proof applies to the space $LM$ 
of {\it continuous} loops and $\widetilde{N}$ still has the structure of a tubular neighborhood and we again have a diffeomorphism $p^*\bigl(D(\nu)\bigr)\cong \widetilde{N}$ between spaces of continuous loops. But we do not need this much here. 

Let $\tilde{u}'=(p\times p)^*(u')\in
\tilde{H}^d(\widetilde{N}/\partial\widetilde{N})$, and $\tilde{u}=(p\times p)^*(u)\in H^d(LM\times LM)$ be pull-backs of Thom classes. Define the
transfer map $j_!$ by the following composition of maps.
\begin{equation} 
j_!: H_*(LM\times LM) \xrightarrow{\tilde{c}_*} 
\tilde{H}_*(\widetilde{N}/\partial\widetilde{N})
\xrightarrow[\cong]{\tilde{u}'\cap(\ )}  H_{*-d}(\widetilde{N}) 
\xrightarrow[\cong]{\tilde{\pi}_*}
H_{*-d}(LM\times_MLM).
\end{equation}
The tubular neighborhood structure of $\widetilde{N}$ implies that the middle map is a genuine Thom isomorphism. 

\begin{definition}\label{definition of loop product} 
Let $M$ be a closed oriented $d$-manifold. For $a,b\in \mathbb H_*(LM)$, their loop product, denoted by $a\cdot b$, is defined by
\begin{equation}\label{loop product} 
a\cdot b=(-1)^{d(|a|'-d)}\iota_*j_!(a\times b)
=(-1)^{d|a|}\iota_*j_!(a\times b). 
\end{equation}
\end{definition} 

The sign $(-1)^{d(|a|'-d)}$ appears in \cite{CJY} in the commutative
diagram (1-7). We include this sign explicitly in the definition of
the loop product for at least three reasons. The most trivial reason
is that on the left hand side, the dot representing the loop product
is between $a$ and $b$. On the right hand side, $j_!$ of degree $-d$
representing the loop product is in front of $a$. Switching $a$ and
$j_!$ gives the sign $(-1)^{d|a|'}$. The other part of the sign
$(-1)^d$ comes from our choice of orientation of $\nu$ and ensures
that $s_*([M])\in \mathbb H_0(LM)$, with the $+$ sign, is the unit of the loop
product. We quickly verify the correctness of the sign.
\begin{lemma} 
The element $s_*([M])\in\mathbb H_0(LM)$ is the unit of the loop
product. Namely for any $a\in \mathbb H_*(LM)$,
\begin{equation}
s_*([M])\cdot a=a\cdot s_*([M])=a.
\end{equation}
\end{lemma}
\begin{proof} We consider the following diagram. 
\begin{equation*}
\begin{CD}
@. LM @= LM  \\
@. @A{\pi_2}AA @| @. \\
M\times M @<{1\times p}<< M\times LM  @<{j'}<< M\times_M LM @= LM  \\
@| @V{s\times 1}VV @V{s\times_M1}VV @|  \\
M\times M  @<{p\times p}<< LM\times LM  @<{j}<< LM\times_M LM @>{\iota}>> LM 
\end{CD}
\end{equation*} 
In the induced homology diagram with transfers $j_!$ and $j'_!$, the
bottom middle square commutes because transfers are defined using Thom
classes pulled back from the same Thom class $u$ of the base
manifold. Thus,
\begin{equation*}
s_*([M])\cdot a=\iota_*j_!(s\times 1)_*([M]\times a)=j'_!([M]\times a). 
\end{equation*}
Here, since $[M]$ has degree $d$, the sign in \eqref{loop product} is $+1$. Since $\pi_2\circ j'=1$, the identity on $LM$,
\begin{equation*}
j'_!([M]\times a)={\pi_2}_*j'_*j'_!([M]\times a)
={\pi_2}_*\bigl((1\times p)^*(u)\cap([M]\times a)\bigr).
\end{equation*}
Due to the way $\nu$ is oriented, the Thom class $u$ is of the form
$u=\{M\}\times 1+\dotsb+(-1)^d1\times\{M\}$. Hence
${\pi_2}_*\bigl((1\times p)^*(u)\cap([M]\times
a)\bigr)={\pi_2}_*([x_0]\times a+\dotsm)=a$.

The other identity $a\cdot s_*([M])=a$ can be proved similarly. This
completes the proof.
\end{proof}

The second reason is that this choice of sign for the loop product is
the same sign appearing in the formula for the intersection product
defined in terms of the transfer map \eqref{intersection
product}. This makes the loop product compatible with the intersection
product in the following sense.

\begin{proposition}
Both of the following maps are algebra maps preserving units between
the loop algebra $\mathbb H_*(LM)$ and the intersection ring $\mathbb H_*(M)$.
\begin{equation}
p_*: \mathbb H_*(LM) \longrightarrow \mathbb H_*(M),\qquad 
s_*: \mathbb H_*(M) \longrightarrow
\mathbb H_*(LM).
\end{equation}
\end{proposition} 
\begin{proof} The proof is more or less straightforward, 
but we discuss it briefly. We consider the following diagram.
\begin{equation*}
\begin{CD}
LM\times LM @<{j}<< LM\times_M LM @>{\iota}>> LM \\ @V{p\times p}VV
@V{p}VV @V{p}VV \\ M\times M @<{\phi}<< M @= M \\ @V{s\times s}VV
@V{s}VV @V{s}VV \\ LM\times LM @<{j}<< LM\times_M LM @>{\iota}>> LM
\end{CD}\end{equation*}
Since the Thom classes for embeddings $j$ and $\phi$ are compatible
via $(p\times p)^*$, the induced homology diagram with transfers
$j_!$ and $\phi_!$ is commutative. Then by diagram chasing, we can
easily check that $p_*$ and $s_*$ preserve products because of the
same signs appearing in \eqref{intersection product} and 
\eqref{loop product}. 
\end{proof} 

The third reason of the sign for the loop product is that it gives the
correct graded commutativity, as given in \cite{CS} proved in terms of chains. We discuss a homotopy theoretic proof of graded commutativity because \cite{CJ} did not include it, and because the homotopy theoretic proof itself is not so trivial with careful treatment of transfers and signs. Contrast the present homotopy theoretic proof with the simple geometric proof given in \cite{CS}. 

\begin{proposition} 
For $a,b\in \mathbb H_*(LM)$, the following graded commutativity relation holds\textup{:}
\begin{equation} 
a\cdot b=(-1)^{(|a|'-d)(|b|'-d)}b\cdot a=(-1)^{|a||b|}b\cdot a. 
\end{equation}
\end{proposition} 
\begin{proof}
We consider the following commutative diagram, where $R_{\frac12}$ is
the rotation of loops by $\frac12$, that is,
$R_{\frac12}(\gamma)(t)=\gamma(t+\frac12)$.
\begin{equation*}
\begin{CD}
LM\times LM @<{j}<< LM\times_M LM @>{\iota}>> LM \\
@V{T}VV    @V{T}VV  @V{R_{\frac12}}VV   \\
LM\times LM @<{j}<< LM\times_M LM @>{\iota}>> LM
\end{CD}
\end{equation*}
Since $R_{\frac12}$ is homotopic to the identity map, we have
${R_{\frac12}}_*=1$. Hence
\begin{equation*}
a\cdot b=(-1)^{d(|a|'-d)}\iota_*j_!(a\times b)
=(-1)^{d(|a|'-d)}\iota_*T_*j_!(a\times b).
\end{equation*} 
Next we show that the induced homology square with transfer $j_!$, we
have $T_*j_!=(-1)^dj_!T_*$. Since the left square in the above diagram
commutes on space level, we have that $T_*j_!$ and $j_!T_*$ coincides
up to a sign. To determine this sign, we compose $j_*$ on the left of
these maps and compare. Since the homology square with induced
homology maps commute,
\begin{equation*}
j_*T_*j_!(a\times b)=T_*j_*j_!(a\times
b)=T_*\bigl(\tilde{u}\cap(a\times b)\bigr).
\end{equation*}
On the other hand, 
\begin{equation*}
j_*j_!T_*(a\times b)=\tilde{u}\cap T_*(a\times b)
=T_*\bigl(T^*(\tilde{u})\cap(a\times b)\bigr). 
\end{equation*}
We compare $T^*(\tilde{u})$ and $\tilde{u}$. Since $\tilde{u}=(p\times
p)^*(u)$, we have $T^*(\tilde{u})=(p\times p)^*T^*(u)$. Since $u$ is
characterized by the property $u\cap[M\times M]=\phi_*([M])$ and
$T\circ \phi=\phi$, we have
\begin{equation*}
\phi_*([M])=T_*\phi_*([M])=T^*(u)\cap T_*([M\times
M])=T^*(u)\cap(-1)^d[M\times M].
\end{equation*}
Thus $T^*(u)=(-1)^du$. Hence $T^*(\tilde{u})=(-1)^d\tilde{u}$. In view
of the above two identities, this implies that
$j_*T_*j_!=(-1)^dj_*j_!T_*$, or $T_*j_!=(-1)^dj_!T_*$. 

Continuing our computation,
\begin{equation*} 
a\cdot b=(-1)^{d|a|'}\iota_*j_!T_*(a\times b)
=(-1)^{|a|'|b|'+d|\alpha|}\iota_*j_!(b\times a)=(-1)^{(|a|'-d)(|b|'-d)}b\cdot a.
\end{equation*}
This completes the homotopy theoretic proof of commutativity formula. 
\end{proof}

\begin{remark}
If we let $\mu=\iota_*j_!:\mathbb H_*(LM)\otimes \mathbb H_*(LM) \longrightarrow \mathbb H_*(LM)$, then using the method in \cite{T2}, we can show that the associativity of $\mu$ takes the form $\mu\circ(1\otimes \mu)=(-1)^d\mu\circ(\mu\otimes 1)$. With our choice of the sign for the loop product in Definition \ref{definition of loop product}, we can get rid of the above sign and we have a usual associativity relation $(a\cdot b)\cdot c=a\cdot(b\cdot c)$ for the loop product without any signs for $a,b,c\in\mathbb H_*(LM)$. This is yet another reason of our choice of the sign in the definition of the loop product. 
\end{remark}

The transfer map $j_!$ enjoys the following properties similar to
those satisfies by $\phi_!$ as given in Proposition~\ref{properties of
phi}. The proof is similar, and we omit it.

\begin{proposition}
For $a,b\in \mathbb H_*(LM)$ and $\alpha\in H^*(LM\times LM)$, the following
formulas are valid.
\begin{align}
j_*j_!(a\times b)&=\tilde{u}\cap(a\times b) \label{j_!1}\\
j_!\bigl(\alpha\cap(a\times b)\bigr)&=
(-1)^{d|\alpha|}j^*(\alpha)\cap j_!(b\times c) \label{j_!2}
\end{align}
\end{proposition}

The second formula says that $j_!$ is a $H^*(LM\times LM)$-module map.

\section{Cap products and extended BV algebra structure}

We examine compatibility of the cap product with the various
structures in the BV-algebra $\mathbb{H}_*(LM)=H_{*+d}(LM)$.

We recall that a BV-algebra $A_*$ is an associative graded commutative
algebra equipped with a degree $1$ Lie bracket $\{\ ,\ \}$ and a
degree $1$ operator $\Delta$ satisfying the following relations for
$a,b,c\in A_*$:
\begin{gather}
\Delta(a\cdot b)=(\Delta a)\cdot b + (-1)^{|a|}a\cdot \Delta b +
(-1)^{|a|}\{a,b\} 
\tag{BV identity} \\ 
\{a, b\cdot c\} =\{a,b\}\cdot c + (-1)^{|b|(|a|+1)}b\cdot\{a,c\} 
\tag{Poisson identity} \\ 
a\cdot b=(-1)^{|a||b|}b\cdot a,\qquad 
\{a,b\}=-(-1)^{(|a|+1)(|b|+1)}\{b,a\}
\tag{Commutativity} \\
\{a,\{b,c\}\}=\{\{a,b\},c\}+(-1)^{(|a|+1)(|b|+1)}\{b,\{a,c\}\}
\tag{Jacobi identity}
\end{gather} 
Here, degrees of elements are given by $\Delta a\in A_{|a|+1}, a\cdot b\in A_{|a|+|b|}$, and $\{a,b\}\in A_{|a|+|b|+1}$. 
One way to view these relations is to consider operators $D_a$ and
$M_a$ acting on $A_*$ for each $a\in A_*$ given by $D_a(b)=\{a,b\}$ and
$M_a(b)=a\cdot b$. Let $[x,y]=xy-(-1)^{|x||y|}yx$ be the graded
commutator of operators. Then the Poisson identity and the Jacobi
identity take the following forms:
\begin{equation}
[D_a,M_b]=M_{\{a,b\}},\qquad [D_a,D_b]=D_{\{a,b\}},
\end{equation}
where degrees of operators are $|D_a|=|a|+1$ and $|M_b|=|b|$. 

One nice context to understand BV identity is in the context of odd symplectic geometry (\cite{G}, \S2), where BV operator $\Delta$ appears as a mixed second order odd differential operator, and BV identity can be simply understood as Leipnitz rule in differential calculus. This context actually arises in loop homology. In \cite{T1}, we explicitly computed the BV structure of $\mathbb H_*(LM)$ for the Lie group $\text{SU}(n+1)$ and complex Stiefel manifolds. There, the BV operator $\Delta$ is given by second order mixed odd differential operator as above, and $\mathbb H_*(LM)$ is interpreted as the space of polynomial functions on the odd symplectic vector space. 

The fact that the loop algebra $\mathbb{H}_*(LM)$ is a
BV-algebra was proved in \cite{CS}. Note that the above BV relations
are satisfied with respect to $\mathbb{H}_*$-grading, rather than the
usual homology grading. 
The same is true for compatibility relations
with cap products. 

First we discuss cohomological $S^1$ action operator $\Delta$ on $H^*(LM)$. Let $\Delta:
S^1\times LM \longrightarrow LM$ be the $S^1$ action map given by
$\Delta(t,\gamma)=\gamma_t$, where $\gamma_t(s)=\gamma(s+t)$ for
$s,t\in S^1=\mathbb{R}/\mathbb{Z}$. The degree $-1$ operator
$\Delta: H^*(LM) \longrightarrow H^{*-1}(LM)$ is defined by the
following formula for $\alpha\in H^*(LM)$:
\begin{equation}
\Delta^*(\alpha)=1\times\alpha + \{S^1\}\times \Delta\alpha
\end{equation}
where $\{S^1\}$ is the fundamental cohomology class of $S^1$. The
homological $S^1$ action $\Delta$ is not a derivation with respect to the loop product and the deviation from being a derivation is given 
by the loop bracket. However, the
cohomology $S^1$-operator $\Delta$ is a derivation with respect to the cup product.

\begin{proposition}
The cohomology $S^1$-operator $\Delta$ satisfies $\Delta^2=0$, and it acts as a derivation on the cohomology ring $H^*(LM)$. That is, for
$\alpha,\beta\in H^*(LM)$,
\begin{equation}\label{delta and cup}
\Delta(\alpha\cup\beta)=(\Delta\alpha)\cup\beta +
(-1)^{|\alpha|}\alpha\cup\Delta\beta.
\end{equation}
\end{proposition}
\begin{proof} The property $\Delta^2=0$ is straightforward 
using the following diagram
\begin{equation*}
\begin{CD}
S^1\times S^1 \times LM @>{1\times\Delta}>> S^1\times LM \\ @V{\mu\times 1}VV
@V{\Delta}VV \\ S^1\times LM @>{\Delta}>> LM
\end{CD}
\end{equation*}
Comparing both sides of $(1\times\Delta)^*\Delta^*(\alpha)=(\mu\times
1)^*\Delta^*(\alpha)$, we obtain $\Delta^2(\alpha)=0$.

For the derivation property, we consider the following diagram. 
\begin{equation*}
\begin{CD}
S^1\times LM @>{\phi\times\phi}>>  (S^1\times S^1)\times (LM\times LM) 
@>{1\times T\times 1}>>  (S^1\times LM)\times (S^1\times LM)  \\
@V{\Delta}VV    @.   @V{\Delta\times\Delta}VV    \\
LM @>{\phi}>> LM\times LM @= LM\times LM 
\end{CD}
\end{equation*} 
On the one hand, $\Delta^*\phi^*(\alpha\times
\beta)=\Delta^*(\alpha\cup\beta) =1\times(\alpha\cup\beta) +
\{S^1\}\times \Delta(\alpha\cup\beta)$. On the other hand,
\begin{equation*}
(\phi\times\phi)^*(1\times T\times 1)^*(\Delta\times
\Delta)^*(\alpha\times \beta) =1\times(\alpha\cup\beta) +
(-1)^{|\alpha|}\{S^1\}\times \bigl(\alpha\cup\Delta\beta 
+ \Delta\alpha\cup\beta\bigr).
\end{equation*}
Comparing the above two identities, we obtain the derivation formula.
\end{proof} 

We can regard the cohomology ring $H^*(LM)$ together with cohomological $S^1$ action $\Delta$ as a BV algebra with trivial bracket product. 

Now we show that the cap product is compatible with the loop product
in the BV-algebra $\mathbb H_*(LM)$. The following theorem describes the behavior of the cap product with those elements in the subalgebra of $H^*(LM)$ generated by $H^*(M)$ and $\Delta\bigl(H^*(M)\bigr)$.

\begin{theorem} Let $\alpha\in H^*(M)$ and $b,c\in \mathbb H_*(LM)$. 
The cap product with $p^*(\alpha)$ behaves
associatively and graded commutatively with respect to the loop
product. Namely
\begin{equation}\label{cap and loop product} 
p^*(\alpha)\cap(b\cdot c)=(p^*(\alpha)\cap b)\cdot c
=(-1)^{|\alpha||b|}b\cdot(p^*(\alpha)\cap c).
\end{equation}

The cap product with 
$\Delta\bigl(p^*(\alpha)\bigr)$ is a derivation on the loop product. Namely, 
\begin{equation}\label{cap derivation} 
\Delta\bigl(p^*(\alpha)\bigr)\cap(b\cdot c)
=\bigl(\Delta\bigl(p^*(\alpha)\bigr)\cap b\bigr)\cdot c +
(-1)^{(|\alpha|-1)|b|}b\cdot\bigl(\Delta\bigl(p^*(\alpha)\bigr)\cap
c\bigr).
\end{equation}
\end{theorem} 
\begin{proof}
For the first formula, we consider the following diagram, 
where $\pi_i$ is the projection onto the $i$th factor for $i=1,2$.  
\begin{equation*}
\begin{CD} 
LM @<{\pi_i}<<  LM\times LM @<{j}<< LM\times_M LM  @>{\iota}>> LM  \\
@V{p}VV  @V{p\times p}VV  @V{q}VV   @V{p}VV \\
M @<{\pi_i}<< M\times M  @<{\phi}<< M @= M 
\end{CD}
\end{equation*} 
Since $p^*(\alpha)\cap(b\cdot c)
=(-1)^{d|b|}\iota_*\bigl(\iota^*p^*(\alpha)\cap j_!(b\times c)\bigr)$, 
we need to understand $\iota^*p^*(\alpha)$. From the above commutative 
diagram, we have $\iota^*p^*(\alpha)=j^*\pi_i^*p^*(\alpha)$, which is
equal to either $j^*(p^*(\alpha)\times 1)$ or $j^*(1\times
p^*(\alpha))$. In the first case, continuting our computation using 
\eqref{j_!2}, we have 
\begin{align*}
p^*(\alpha)\cap(b\cdot c)&=
(-1)^{d|b|}\iota_*\bigl(j^*(p^*(\alpha)\times 1)
\cap j_!(b\times c)\bigr) \\
&=(-1)^{d|b|+d|\alpha|}\iota_*j_!
\bigl((p^*(\alpha)\times 1)\cap(b\times c)\bigr) \\
&=(-1)^{d|b|+d|\alpha|}\iota_*j_!
\bigl((p^*(\alpha)\cap b)\times c\bigr)  \\
&=(p^*(\alpha)\cap b)\cdot c.
\end{align*}
Similarly, using $\iota^*p^*(\alpha)
=j^*\bigl(1\times p^*(\alpha)\bigr)$, we get the other identity. This
proves (1). 

For (2), we first note that the element
$\Delta\bigl(p^*(\alpha)\bigr)\cap(b\cdot c)$ is equal to 
\begin{equation*}
\Delta\bigl(p^*(\alpha)\bigr)\cap (-1)^{d|b|}\iota_*j_!(b\times c)
=(-1)^{d|b|}\iota_*\bigl(\iota^*\bigl(\Delta(p^*(\alpha))\bigr)
\cap j_!(b\times c)\bigr).
\end{equation*}
Thus, we need to understand the element
$\iota^*\bigl(\Delta(p^*(\alpha))\bigr)$. We need some notations. Let
$I=I_1\cup I_2$, where $I_1=[0,\frac12]$ and $I_2=[\frac12,1]$, and
set $S_i^1=I_i/\partial I_i$ for $i=1,2$. Let $r: S^1=I/\partial I 
\longrightarrow I/\{0,\frac12,1\}=S^1_1\vee S^1_2$ be an
identification map, and let $\iota_i: S^1_i \longrightarrow S^1_1\vee
S^1_2$ be the inclusion map into the $i$th wedge summand. We examine
the following diagram 
\begin{equation*}
\begin{CD}
S^1\times(LM\underset{M}{\times}LM)  @>{r\times 1}>> 
(S^1_1\vee S^1_2)\times
(LM\underset{M}{\times}LM) @<<< \{0\}\times (LM\underset{M}{\times}LM)  \\
@V{1\times\iota}VV   @V{e'}VV   @V{\iota}VV   \\
S^1\times LM  @>{e}>> M @<{p}<<  LM 
\end{CD}
\end{equation*} 
where $e=p\circ\Delta$ is the evaluation map for $S^1\times LM$, and
the other evaluation map $e'$ is given by 
\begin{equation*}
e'(t,\gamma,\eta)=
\begin{cases}
\gamma(2t)& 0\le t\le\frac12,\\
\eta(2t-1)& \frac12\le t\le 1.
\end{cases}
\end{equation*}
For $\alpha\in H^*(M)$, we let 
\begin{equation*}
{e'}^*(\alpha)=1\times \iota^*p^*(\alpha) 
+\{s^1_1\}\times \Delta_1(\alpha) 
+\{S^1_2\}\times \Delta_2(\alpha)
\end{equation*}
for some $\Delta_i(\alpha)\in H^*(LM\times_M LM)$ for $i=1,2$. The
first term in the right hand side is identified using the right square
of the above commutative diagram. Since $r^*(\{S^1_i\})=\{S^1\}$ for
$i=1,2$,  
\begin{equation*}
(r\times 1)^*{e'}^*(\alpha)=1\times \iota^*p^*(\alpha) + 
\{S^1\}\times \bigl(\Delta_1(\alpha)+\Delta_2(\alpha)\bigr). 
\end{equation*}
The commutativity of the left square implies that this must be equal to 
\begin{equation*}
(1\times\iota)^*\Delta^*p^*(\alpha)
=1\times \iota^*p^*(\alpha)+\{S^1\}\times 
\iota^*\Delta\bigl(p^*(\alpha)\bigr). 
\end{equation*}
Hence we have 
\begin{equation*}
\iota^*\bigl(\Delta(p^*(\alpha))\bigr)
=\Delta_1(\alpha)+\Delta_2(\alpha)\in H_*(LM\times_M LM). 
\end{equation*}
To understand elements $\Delta_i(\alpha)$, we consider the following
commutative diagram, where $\ell_1(t)=2t$ for $0\le t\le\frac12$ and 
$\ell_2(t)=2t-1$ for $\frac12\le t\le 1$. 
\begin{equation*}
\begin{CD}
S^1_i\times(LM\times_M LM) @>{\ell_i\times j}>> S^1\times(LM\times LM)
@>{1\times\pi_i}>> S^1\times LM \\
@V{\iota_i\times 1}VV @. @V{\Delta}VV \\
(S^1_1\vee S^1_2)\times(LM\times_M LM) @>{e'}>> M  @<{p}<<  LM 
\end{CD} 
\end{equation*} 
On the one hand, $(\iota_1\times 1)^*{e'}^*(\alpha)
=1\times\iota^*p^*(\alpha)+ \{S^1_1\}\times\Delta_1(\alpha)$. On the
other hand, 
\begin{equation*}
(\ell_1\times j)^*(1\times \pi_1)^*\Delta^*p^*(\alpha)=
1\times j^*\bigl(p^*(\alpha)\times 1\bigr)
+\{S^1_1\}\times j^*\bigl(\Delta(p^*(\alpha))\times 1\bigr).
\end{equation*}
By the commutativity of the diagram, we
get $\Delta_1(\alpha)=j^*\bigl(\Delta(p^*(\alpha))\times
1\bigr)$. Similarly, $i=2$ case implies $\Delta_2(\alpha)
=j^*\bigl(1\times\Delta(p^*(\alpha))\bigr)$. Hence we finally obtain 
\begin{equation*}
\iota^*\bigl(\Delta(p^*(\alpha))\bigr)
=j^*\bigl(\Delta(p^*(\alpha))\times 1 + 1\times
\Delta(p^*(\alpha))\bigr).
\end{equation*}
With this identification of $\iota^*\bigl(\Delta(p^*(\alpha))\bigr)$
as $j^*$ of some other element, we can continue our initial
computation. 
\begin{equation*}
\begin{split}
\Delta\bigl(&p^*(\alpha)\bigr)\cap(b\cdot c)
=(-1)^{d|b|}\iota_*\bigl(j^*\bigl(\Delta(p^*(\alpha))\times 1 + 
1\times \Delta\bigl(p^*(\alpha)\bigr)\bigr)\cap j_!(b\times c)\bigr)
\\
&=(-1)^{d|b|+(|\alpha|-1)d} 
\iota_*j_!\Bigl(\bigl(\Delta(p^*(\alpha))\times 1 
+1\times \Delta\bigl(p^*(\alpha)\bigr)\bigr)\cap(b\times c)\Bigr) \\
&=(-1)^{d(|\alpha|+|b|-1)}\iota_*j_!\Bigl(
\bigl(\Delta(p^*(\alpha))\cap b\bigr)\times c +
(-1)^{(|b|+d)(|\alpha|-1)}b\times 
\bigl(\Delta(p^*(\alpha))\cap c\bigr)\Bigr) \\
&=\bigl(\Delta(p^*(\alpha))\cap b\bigr)\cdot c +
(-1)^{(|\alpha|-1)|b|}b\cdot\bigl(\Delta(p^*(\alpha))\cap c\bigr).
\end{split}
\end{equation*}
This completes the proof of the derivation property of the cap product
with respect to the loop product. 
\end{proof} 

Next we describe the relation between the cap product and the BV
operator in homology and cohomology. 

\begin{proposition} For $\alpha\in H^*(LM)$ and $b\in \mathbb H_*(LM)$, the BV-operator $\Delta$ satisfies 
\begin{equation}\label{delta and cap}
\Delta(\alpha\cap b)=(\Delta\alpha)\cap b+(-1)^{|\alpha|}\alpha\cap
\Delta b.
\end{equation}
\end{proposition}
\begin{proof}
On the one hand, the $S^1$-action map $\Delta: S^1\times LM
\longrightarrow LM$ satisfies 
\begin{equation*}
\Delta_*\bigl(\Delta^*(\alpha)\cap([S^1]\times b)\bigr)
=\alpha\cap\Delta_*([S^1]\times b)
=\alpha\cap\Delta b.
\end{equation*}
On the other hand, since $\Delta^*(\alpha)=1\times\alpha +
\{S^1\}\times \Delta\alpha$, we have 
\begin{equation*}
\begin{split}
\Delta_*\bigl(\Delta^*(\alpha)\cap([S^1]\times b)\bigr)
&=\Delta_*\bigl(
(-1)^{|\alpha|}[S^1]\times(\alpha\cap b) +
(-1)^{|\alpha|-1}[pt]\times(\Delta\alpha\cap b)\bigr)  \\
&=(-1)^{|\alpha|}\Delta(\alpha\cap b) + 
(-1)^{|\alpha|-1}\Delta\alpha\cap b. 
\end{split}
\end{equation*}
Comparing the above two formulas, we obtain 
$\Delta(\alpha\cap b)=\Delta\alpha\cap b + 
(-1)^{|\alpha|}\alpha\cap\Delta b$. 
\end{proof} 
Since homology BV operator $\Delta$ on $\mathbb H_*(LM)$ acts trivially on $\mathbb H_*(M)$, the
following corollary is immediate. 
\begin{corollary}
For $\alpha\in H^*(M)$, the cap product of $\Delta\alpha$ with
$\mathbb H_*(M)\subset \mathbb H_*(LM)$ is trivial.
\end{corollary} 
\begin{proof} For $b\in\mathbb H_*(M)$, the operator $\Delta$ acts trivially on both $\alpha\cap b$ and $b$. Hence formula \eqref{delta and cap} implies $(\Delta\alpha)\cap b=0$. 
\end{proof} 

Next, we discuss a behavior of the cap product with respect to the
loop bracket.

\begin{theorem}  The cap product with $\Delta\bigl(p^*(\alpha)\bigr)$ 
is a derivation on the loop bracket. Namely, for $\alpha\in H^*(M)$
and $b,c\in \mathbb H_*(LM)$,
\begin{equation}\label{cap-loop bracket}
\Delta\bigl(p^*(\alpha)\bigr)\cap\{b,c\}
=\{\Delta\bigl(p^*(\alpha)\bigr)\cap b, c\} +
(-1)^{(|\alpha|-1)(|b|-1)}\{b, \Delta\bigl(p^*(\alpha)\bigr)\cap
c\}.
\end{equation}
\end{theorem} 
\begin{proof}  Our proof is computational using previous results. 
We use the BV identity as the definition of the loop bracket. Thus,
\begin{equation*}
\{b,c\}=(-1)^{|b|}\Delta(b\cdot c)-(-1)^{|b|}(\Delta b)\cdot c 
-b\cdot \Delta c.
\end{equation*} 
We compute the right hand side of \eqref{cap-loop bracket}. 
For simplicity, we write $\Delta\alpha$
for $\Delta\bigl(p^*(\alpha)\bigr)$.  Each term in the right hand side
of \eqref{cap-loop bracket} gives
\begin{gather*}
\!\!\!\!\!\!\!\{\Delta\alpha\cap b, c\}
=(-1)^{|b|-|\alpha|+1}\Delta\bigl((\Delta\alpha\cap b)\cdot c\bigr)
-(-1)^{|b|}(\Delta\alpha\cap \Delta b)\cdot c 
-(\Delta\alpha\cap b)\cdot\Delta c, \\
\!\!\!\!\!\!\!\!\!\!\!\!\{b,\Delta\alpha\cap c\}
=(-1)^{|b|}\Delta\bigl(b\cdot(\Delta\alpha\cap c)\bigr)
-(-1)^{|b|}\Delta b\cdot(\Delta\alpha\cap c) 
-(-1)^{|\alpha|-1}b\cdot(\Delta\alpha\cap\Delta c),
\end{gather*}
Here we used \eqref{delta and cap} for the second term in the first
identity and in the third term in the second identity. Combining these
formulas, we get
\begin{multline*} 
\{\Delta\alpha\cap b, c\}+ (-1)^{(|\alpha|-1)(|b|+1)}
\{b,\Delta\alpha\cap c\}  \\
=\bigl[ 
(-1)^{|b|-|\alpha|+1}\Delta\bigl((\Delta\alpha\cap b)\cdot c\bigr)
+(-1)^{|b|+(|\alpha|-1)(|b|+1)}
\Delta\bigl(b\cdot(\Delta\alpha\cap c)\bigr)\bigr]  \\
-\bigl[(-1)^{|b|}(\Delta\alpha\cap\Delta b)\cdot c + 
(-1)^{(|\alpha|-1)(|b|+1)+|b|}\Delta b\cdot(\Delta\alpha\cap c)\bigr] \\
-\bigl[(\Delta\alpha\cap b)\cdot \Delta c +
(-1)^{(|\alpha|-1)(|b|+1)+|\alpha|-1}b\cdot(\Delta\alpha\cap \Delta c)\bigr].
\end{multline*}
Using the derivation formula for $\Delta\alpha\cap(\ )$ with respect
to the loop product \eqref{cap derivation}, three pairs of terms above
become
\begin{multline*}
(-1)^{|b|-|\alpha|+1}\Delta\bigl(\Delta\alpha\cap(b\cdot c)\bigr)
-(-1)^{|b|}\Delta\alpha\cap(\Delta b\cdot c)
-\Delta\alpha\cap(b\cdot\Delta c) \\
=\Delta\alpha\cap\bigl[(-1)^{|b|}\Delta(b\cdot c)
-(-1)^{|b|}\Delta b\cdot c -b\cdot\Delta c\bigr] = \Delta\alpha\cap\{b,c\}.
\end{multline*}
This completes the proof of the derivation formula for the loop bracket. 
\end{proof} 

Recall that in the BV algebra $\mathbb H_*(LM)$, for every $a\in \mathbb H_*(LM)$ the
operation $\{a,\ \cdot \ \}$ of taking the loop bracket with $a$ is a
derivation with respect to both the loop product and the loop bracket,
in view of the Poisson identity and the Jacobi identity. Since we have
proved that the cap product with $\Delta p^*(\alpha)$ for $\alpha\in
H^*(M)$ is a derivation with respect to both the loop product and the
loop bracket, we wonder if we can extend the BV structure in $\mathbb H_*(LM)$ to a BV structure in $H^*(M)\oplus \mathbb H_*(LM)$. Indeed this is possible by extending the loop product and the loop bracket to elements in $H^*(M)$ as follows.
\begin{definition}
For $\alpha,\beta\in H^*(M)$ and $b\in \mathbb H_*(LM)$, we define their loop product and loop bracket by
\begin{equation}
\begin{gathered}
\alpha\cdot b=\alpha\cap b,\qquad 
\{\alpha,b\}=(-1)^{|\alpha|}(\Delta\alpha)\cap b,\\
\alpha\cdot\beta=\alpha\cup\beta, \qquad \{\alpha,\beta\}=0.
\end{gathered}
\end{equation}
This defines an associative graded commutative loop product by \eqref{cap and loop product}, and a bracket product on $H^*(M)\oplus\mathbb H_*(LM)$. 
\end{definition}
Note that this loop product on $H^*(M)\oplus\mathbb H_*(LM)$ reduces to the ring structure on $H^*(M)\oplus \mathbb H_*(M)$ mentioned in the introduction. 

With this definition, Poisson identities and Jacobi identities are
still valid in $H^*(M)\oplus \mathbb H_*(LM)$.

\begin{theorem} Let $\alpha,\beta\in H^*(M)$, and let $b,c\in \mathbb H_*(LM)$. 

\noindent\textup{(I)} The following Poisson identities are valid in $H^*(M)\oplus\mathbb H_*(LM)$\textup{:} 
\begin{align}
\{\alpha,\beta\cdot c\}
&=\{\alpha,\beta\}\cdot c+ 
(-1)^{|\beta|(|\alpha|-1)}\beta\cdot\{\alpha,c\} \label{eq1}\\
\{\alpha\beta,c\}&=\alpha\cdot\{\beta,c\}
+(-1)^{|\alpha||\beta|}\beta\cdot\{\alpha,c\} \label{eq2}\\
\{\alpha,b\cdot c\}&=\{\alpha,b\}\cdot c
+(-1)^{(|b|-d)(|\alpha|-1)}b\cdot\{\alpha,c\} \label{eq3}\\
\{\alpha\cdot b,c\}&=\alpha\cdot\{b,c\}+
(-1)^{|\alpha|(|b|-d)}b\cdot\{\alpha,c\}. \label{eq4} 
\end{align}

\noindent\textup{(II)} The following Jacobi identities are valid in $H^*(M)\oplus\mathbb H_*(LM)$\textup{:}
\begin{align}
\{\alpha,\{\beta,c\}\}&=\{\{\alpha,\beta\},c\}+
(-1)^{(|\alpha|-1)(|\beta|-1)}\{\beta,\{\alpha,c\}\} \label{eq5} \\
\{\alpha,\{b,c\}\}&=\{\{\alpha,b\},c\}+
(-1)^{(|\alpha|-1)(|b|-d+1)} \{b,\{\alpha,c\}\}.  \label{eq6}
\end{align}
\end{theorem}
\begin{proof}  If we unravel definitions, we see that \eqref{eq1} 
and \eqref{eq5} are really the same as the graded commutativity of the
cup product of the following form
\begin{align*}
(\Delta\alpha)\cap(b\cap c)&=
(-1)^{|\beta|(|\alpha|-1)}\beta\cap(\Delta\alpha\cap c), \\
(\Delta\alpha)\cap(\Delta\beta\cap c)&=(-1)^{(|\alpha|-1)(|\beta|-1)}
(\Delta\beta)\cap\bigl((\Delta\alpha)\cap c\bigr).
\end{align*}
the identity \eqref{eq2} is equivalent to the derivation formula \eqref{delta and cup} of
the cohomology $S^1$ action operator with respect to the cup product.
\begin{equation*}
\Delta(\alpha\cup \beta)=(\Delta\alpha)\cup \beta + (-1)^{|\alpha|}
\alpha\cup (\Delta \beta).
\end{equation*}
The identity \eqref{eq3} says that $\Delta\alpha\cap (\ )$ is a
derivation with respect to the loop product, and the identity
\eqref{eq6} says that $\Delta\alpha\cap(\ )$ is a derivation with
respect to the loop bracket. We have already verified both of these
cases. Thus, what remains to be checked is formula \eqref{eq4}, which
says
\begin{equation*}
\{\alpha\cap b,c\}=\alpha\cap\{b,c\}+
(-1)^{|\alpha||b|+|\alpha|}b\cdot(\Delta\alpha\cap c).
\end{equation*}
Using the BV identity, the derivation formula \eqref{delta and cap} 
of the BV operator with
respect to the cap product, and properties of $\alpha\cap(\ )$ and
$\Delta\alpha\cap(\ )$, we can prove this identity as follows.
\begin{multline*}
(-1)^{|b|-|\alpha|}\{\alpha\cap b,c\}
=\Delta\bigl((\alpha\cap b)\cdot c\bigr)-\Delta(\alpha\cap b)\cdot c
-(-1)^{|b|-|\alpha|}(\alpha\cap b)\cdot\Delta c \\
=\Delta\bigl(\alpha\cap(b\cdot c)\bigr)
-(\Delta\alpha\cap b+(-1)^{|\alpha|}\alpha\cap\Delta b)\cdot c
-(-1)^{|b|-|\alpha|}\alpha\cap(b\cdot\Delta c) \\
=(\Delta\alpha)\cap(b\cdot c)
-(\Delta\alpha\cap b)\cdot c
+(-1)^{|\alpha|}\alpha\cap\Delta(b\cdot c) \\
-(-1)^{|\alpha|}\alpha\cap(\Delta b\cdot c)
-(-1)^{|b|-|\alpha|}\alpha\cap(b\cdot\Delta c) \\
=(-1)^{(|\alpha|-1)|b|}b\cdot(\Delta\alpha\cap c)
+(-1)^{|\alpha|+|b|}\alpha\cap\{b,c\}.
\end{multline*}
Canceling some signs, we get the desired formula. This completes the proof. 
\end{proof}

Other Poisson and Jacobi identities with cohomology elements in the second argument formally follow from above identities by making following definitions for $\alpha\in H^*(M)$ and $b\in\mathbb H_*(LM)$: 
\begin{equation*}
b\cdot\alpha=(-1)^{|\alpha||b|}\alpha\cdot b,\qquad 
\{b,\alpha\}=-(-1)^{(|\alpha|+1)(|b|+1)}\{\alpha,b\}.
\end{equation*}

For $\alpha\in H^*(M)$ we showed that $\Delta\alpha\cap(\ )$ is a
derivation for both the loop product and the loop bracket, and
$\alpha\cap(\ )$ is graded commutative and associative with respect to
the loop product.  What is the behavior
of $\alpha\cap(\ )$ is with respect to the loop bracket? 
Formula \eqref{eq4} says that $\alpha\cap(\ \cdot\ )$ on loop bracket  is not a derivation or graded commutativity: it is a Poisson identity!

Poisson identities and Jacobi identities we have just proved in $A_*=H^*(M)\oplus\mathbb H_*(LM)$ show that $A_*$ is a Gerstenhaber algebra. In fact, $A_*$ can be formally turned into a BV algebra by defining a BV operator $\boldsymbol\Delta$ on $A_*$ to be trivial on $H^*(M)$ and to be the usual one on $\mathbb H_*(LM)$ coming from the homological $S^1$ action. 

\begin{corollary}\label{BV structure on direct sum} 
The direct sum $A_*=H^*(M)\oplus\mathbb H_*(LM)$ has the structure of a BV algebra. 
\end{corollary}
\begin{proof} Since $\mathbb H_*(LM)$ is a BV algebra and since we have already verified Poisson identities and Jacobi identities in $A_*$, we only have to verify BV identities in $A_*$. For $\alpha,\beta\in H^*(M)$, an identity 
\begin{equation*}
\boldsymbol\Delta(\alpha\cup\beta)=(\boldsymbol\Delta\alpha)\cup\beta
+(-1)^{|\alpha|}\alpha\cup(\boldsymbol\Delta\beta)
+(-1)^{|\alpha|}\{\alpha,\beta\}
\end{equation*}
is trivially satisfied since all terms are zero by definition of BV operator $\boldsymbol\Delta$ and the loop bracket on $H^*(M)\subset A_*$. 

Next, let $\alpha\in H^*(M)$ and $b\in\mathbb H_*(LM)$. Since the BV operator $\boldsymbol\Delta$ on $A_*$ acts trivially on $H^*(M)$, an identity 
\begin{equation*}
\boldsymbol\Delta(\alpha\cap b)=(\boldsymbol\Delta\alpha)\cap b
+(-1)^{|\alpha|}\alpha\cap(\boldsymbol\Delta b)
+(-1)^{|\alpha|}\{\alpha,b\}
\end{equation*}
is really a restatement of the derivative formula of the homology $S^1$ action operator $\Delta$ on cap product: $\Delta(\alpha\cap b)=(-1)^{|\alpha|}\alpha\cap(\Delta b)+(\Delta\alpha)\cap b$ in formula \eqref{delta and cap}. 
\end{proof}

In connection with the above Corollary, we can ask whether $H^*(LM)\oplus \mathbb H_*(LM)$ has a structure of a BV algebra. Of course, $H^*(LM)$ together with the cohomological $S^1$ action operator $\Delta$, which is a derivation, is a BV algebra with trivial bracket product. Thus, as a direct sum of BV algebras, $H^*(LM)\oplus\mathbb H_*(LM)$ is a BV algebra, although products between $H^*(LM)$ and $\mathbb H_*(LM)$ are trivial. More meaningful question would be to ask whether the direct sum $H^*(LM)\oplus\mathbb H_*(LM)$ has a BV algebra structure extending the one on $A_*$ described in Corollary \ref{BV structure on direct sum}. If we want to use the cap product as an extension of the loop product, the answer is no. This is because the cap product with an arbitrary element $\alpha\in H^*(LM)$ does not behave associatively with respect to the loop product in $\mathbb H_*(LM)$: if $\alpha$ is of the form $\alpha=\Delta\beta$ for some $\beta\in H^*(M)$, then $\alpha\cap(\ \cdot\ )$ acts as a derivation on loop product in $\mathbb H_*(LM)$ due to \eqref{cap derivation} and does not satisfy associativity.  

\begin{remark} In the course of our investigation, we noticed the 
following curious identity, which is in some sense symmetric in three
variables, for $\alpha\in H^*(M)$ and $b,c\in\mathbb H_*(LM)$. 
\begin{equation}
\begin{split}
\{\alpha,b\cdot c\}+(-1)^{|b|}\alpha\cdot\{b,c\}
&=\{\alpha,b\}\cdot c+(-1)^{|b|}\{\alpha\cdot b,c\} \\
&=(-1)^{(|\alpha|+1)|b|}\bigl(b\cdot\{\alpha,c\}
+(-1)^{|\alpha|}\{b,\alpha\cdot c\}\bigr).
\end{split}
\end{equation}
This identity is easily proved using Poisson identities. But we wonder
the meaning of this symmetry.
\end{remark}

\section{Cap products in terms of BV algebra structure}

In the previous section, we showed that the BV algebra structure in $\mathbb H_*(LM)$ can be extended to the BV algebra structure in $H^*(M)\oplus \mathbb H_*(LM)$ by proving Poisson identities and Jacobi identities. This may be a bit surprising. But this turns out to be very natural through Poincar\'e duality in the following way. For $a\in \mathbb H_*(M)$, we denote the element $s_*(a)\in \mathbb H_*(LM)$ by $a$, where $s:M\to LM$ is the inclusion map.  

\begin{theorem} For $a\in \mathbb H_*(M)$, let $\alpha=D(a)\in H^*(M)$ be its Poincar\'e dual. Then for any $b\in \mathbb H_*(LM)$, the following identities hold. 
\begin{equation}
p^*(\alpha)\cap b=a\cdot b,\qquad 
(-1)^{|\alpha|}\Delta \bigl(p^*(\alpha)\bigr)\cap b=\{a,b\}.
\end{equation}
\end{theorem}
\begin{proof} Let $1=s_*([M])\in \mathbb H_0(LM)$ be the unit of the loop product. Since $p^*(\alpha)\cap b=p^*(\alpha)\cap(1\cdot b)=
\bigl(p^*(\alpha)\cap1\bigr)\cdot b$ by \eqref{cap and loop product}, and since 
\begin{equation*}
p^*(\alpha)\cap1=p^*(\alpha)\cap s_*([M])=s_*\bigl(s^*p^*(\alpha)\cap[M]\bigr)
=s_*(\alpha\cap[M])=a,
\end{equation*}
we have $p^*(\alpha)\cap b=a\cdot b$. This proves the first identity. 

For the second identity, in the BV identity
\begin{equation*}
(-1)^{|a|}\{a,b\}=\Delta(a\cdot b)-(\Delta a)\cdot b
-(-1)^{|a|}a\cdot\Delta b,
\end{equation*}
the first term in the right hand side gives 
\begin{equation*}
\Delta(a\cdot b)=\Delta(p^*(\alpha)\cap
b)=\Delta\bigl(p^*(\alpha)\bigr)\cap b
+(-1)^{|\alpha|}p^*(\alpha)\cap\Delta b
\end{equation*}
in view of the first identity we just proved and the derivation
property of the homological $A^1$ action operator on cap products. Here
$p^*(\alpha)\cap\Delta b=a\cdot\Delta b$. Since $a\in \mathbb H_*(M)$ is a
homology class of constant loops, we have $\Delta a=0$. Thus,
\begin{equation*}
(-1)^{|a|}\{a,b\}=\Delta\bigl(p^*(\alpha)\bigr)\cap b
+(-1)^{|\alpha|}a\cdot\Delta b-(-1)^{|a|}a\cdot\Delta b
=\Delta\bigl(p^*(\alpha)\bigr)\cap b,
\end{equation*}
since $|\alpha|=-|a|$. Thus, $\{a,b\}=(-1)^{|\alpha|}
\Delta\bigl(p^*(\alpha)\bigr)\cap b$. This completes the proof. 
\end{proof}

In view of this theorem, since $\mathbb H_*(LM)$ is already a BV algebra, the
Poisson identities and Jacobi identities we proved in section 4 may
seem obvious. However, what we did in section 4 is that we gave a
\emph{new and elementary homotopy theoretic proof} of Poisson
identities and Jacobi identities using only basic properties of the
cap product and the BV identity, when at least one of the elements is
from $\mathbb H_*(M)$.

The above theorem shows that loop products and loop brackets with elements in $\mathbb H_*(M)$ can be written as cap products with cohomology elements in $LM$. Thus, compositions of loop products and loop brackets with elements in $\mathbb H_*(M)$ corresponds to a cap product with the product of corresponding cohomology classes in $H^*(LM)$. Namely,
\begin{corollary}
Let $a_0, a_1, \dotsc, a_r\in \mathbb H_*(M)$, and let
$\alpha_0,\alpha_1,\dotsc\alpha_r\in H^*(M)$ be their Poincar\'e
duals. Then for $b\in \mathbb H_*(LM)$,
\begin{equation}
a_0\cdot\{a_1,\{a_2,\dotsc\{a_r,b\}\dotsb\}\}
=(-1)^{|a_1|+\dotsb+|a_r|}
\bigl[\alpha_0(\Delta\alpha_1)(\Delta\alpha_2)
\dotsm(\Delta\alpha_r)\bigr]\cap b.
\end{equation}
\end{corollary}

In section 2, we considered a problem of intersections of loops with
submanifolds in certain configurations, and we saw that the homology
class of the intersections of interest can be given by a cap product
with cohomology cup products of the above form (Proposition \ref{loop
intersection}). The above corollary computes this homology class in
terms of BV structure in $\mathbb H_*(LM)$ using the homology classes of these
submanifolds.

\begin{remark} \label{exterior algebra}
In general, elements $\alpha,\Delta\alpha$ for $\alpha\in H^*(M)$ do not generate the entire cohomology ring $H^*(LM)$. However, if $H^*(M;\mathbb{Q})=\Lambda_{\mathbb{Q}}(\alpha_1,\alpha_2,\dotsc\alpha_r)$ is an exterior algebra, over $\mathbb{Q}$, then using minimal models or spectral sequences, we have
\begin{equation}
H^*(LM;\mathbb{Q})=\Lambda_{\mathbb{Q}}(\alpha_1,\alpha_2,\dotsc\alpha_r)\otimes
\mathbb{Q}[\Delta\alpha_1,\Delta\alpha_2,\dotsc\Delta\alpha_r],
\end{equation}
and thus we have the complete description of the cap products with any elements in $H^*(LM;\mathbb{Q})$ in terms of the BV structure in $\mathbb H_*(LM;\mathbb{Q})$. 
\end{remark}

\end{document}